\documentclass[12pt,a4paper,reqno]{amsart}

\usepackage{amssymb}
\usepackage{amsthm}

\usepackage{amscd}
\usepackage{amsfonts}

\usepackage{amsmath}
\usepackage{amsfonts}
\usepackage{latexsym}
\usepackage{graphicx}
\usepackage{enumerate}
\usepackage{color}
\usepackage{comment}

\usepackage[english]{babel}
\usepackage[latin1]{inputenc}
\newcommand{\N}{\ensuremath{\mathbb{N}}}
\newcommand{\Z}{\ensuremath{\mathbb{Z}}}

\newcommand{\C}{\ensuremath{\mathbb{C}}}

\def\ds{\displaystyle}
\newcommand{\rme}{{\mathrm{e}}}

\newcommand{\eps}{\varepsilon}

\newcommand{\VE}{\mathrm{VE}}
\newcommand{\LVE}{\mathrm{LVE}}
\newcommand{\vfi}{\varphi}

\input xy
\xyoption{all}
\def\sympow{{\setbox0\hbox{$\bigcirc$}\setbox1\hbox to\wd0{\hss$s$\hss}%
\wd1 0pt\box1\box0}}%symmetric power

\newtheorem{theorem}{Theorem}[section]

\newtheorem{proposition}[theorem]{Proposition}
\newtheorem{cor}[theorem]{Corollary}
\newtheorem{remark}[theorem]{Remark}

\newtheorem{definition}[theorem]{Definition}

\def\ds{\displaystyle}
\def\rme{\rm e}
\def\proof{\noindent {\sf Proof.}  }
\def\qed{\hfill $\Box$ \\ \bigskip}

\title[Differential Galois Theory \& Planar Polynomial Vector Fields]{Differential Galois Theory and non-Integrability of Planar Polynomial Vector Fields}

\author[P.~Acosta-Hum\'anez]{Primitivo B. Acosta-Hum\'anez}
\address[P.~Acosta-Hum\'anez]{Facultad de Ciencias B\'asicas y Biom\'edicas, Universidad Sim\'on Bol\'{\i}var, Barranquilla - Colombia}
\email{primitivo.acosta@unisimonbolivar.edu.co -- primi@intelectual.co}

\author[J.~T. L\'azaro]{J.~Tom\'as L\'azaro}
\address[J.~T. L\'azaro]{Departament de Matem\`atiques, Universitat Polit\`ecnica de Catalunya, Spain.}
\email{jose.tomas.lazaro@upc.edu}

\author[J.J.~Morales-Ruiz]{Juan J. Morales-Ruiz}
\address[J.J.~Morales-Ruiz]{Departamento de Matemática Aplicada, Universidad Politécnica de Madrid, Spain.}
\email{juan.morales-ruiz@upm.es}

\author[Ch.~Pantazi]{Chara Pantazi}
\address[Ch.~Pantazi]{Departament de Matem\`atiques, Universitat Polit\`ecnica de Catalunya, Spain.}
\email{chara.pantazi@upc.edu}


\begin{document}

\begin{abstract}
We study a necessary condition for the integrability of the polynomials fields in the plane by means of the differential Galois theory. More concretely, by means of the variational equations  around a particular solution it is obtained a   necessary condition for the existence of a rational  first integral. The method is systematic starting with the first order variational equation.  We illustrate this result with several families of examples. A key point  is to check wether a suitable primitive is elementary or not. Using a theorem by Liouville, the problem is equivalent to the existence of a rational solution of a certain first order linear equation, the Risch equation. This is a classical  problem  studied by Risch in 1969, and the solution is given by the ``Risch algorithm". In this way we point out the connection of the non integrablity with some higher transcendent functions, like the error function.
\end{abstract}

\date{\today}

\maketitle

\section{Introduction}

The problem of the \emph{integrability} of planar vector fields has attracted the attention of many mathematicians during decades. Among the different approaches, the \emph{Galois theory} of linear differential equations has played an important r\^{o}le in its understanding, even in the \emph{a priori} (so far) simpler case of polynomial vector fields (see~\cite{CH,SV,almp} and references therein).
For instance, the application of differential Galois theory to variational equations along a given integral curve constitutes a powerful criterium of non-integrability for Hamiltonian systems (see~\cite{morales 1}).
Ayoul and Zung~\cite {AZ} extended this method to the study of some kind of non Hamiltonian fields. They strongly relied on the main result of~\cite{MRS}.

The aim of this paper is to apply these ideas to prove the non-existence of rational first integrals for planar polynomial vector fields.
%In a previous paper (\cite{almp}) we studied the integrability of some fields on the plane.
Let us consider planar vector fields of the form
\begin{equation}\label{campo}
X=P\dfrac{\partial }{\partial x}+Q \dfrac{\partial }{\partial y},
\end{equation}
with $P,Q$ analytic functions in some domain of $\C^2$ and assume $\Gamma: y-\varphi(x)=0$ to be an integral curve\footnote{Usually it is also refered as an \emph{orbit} of the ode system $\dot{x}=P, \dot{y}=Q$. In qualitative theory of dynamical systems this is commonly called \emph{invariant} curve.} of $X$. This is equivalent to say that
$\Gamma$ is a solution (a \emph{leaf}) of the first order differential equation
\begin{equation} \label{foliacion}
y^\prime=\frac QP=f(x,y),
\end{equation}
which defines its associated foliation (orbits of the vector field $X$). From now on $\phantom{x}'$ will denote derivative with respect to the spatial variable $x$.
The behaviour around the solution $\Gamma$ is usually approached by studying its variational equations. In our case, with respect to equation~\eqref{foliacion}. Precisely,
let's $\phi(x,y)$ denote the flow of~\eqref{foliacion}. Consider $(x_0,y_0)$ a point in $\Gamma$, that is $y_0=\varphi(x_0)$.
Note that $(x_0,y_0)$ is the initial condition defining $\Gamma$: $\varphi(x)=\phi(x,y_0) $. We are interested in the variation of the flow $\phi$
respect to the initial condition $y$, around $y=y_0$ and keeping $x=x_0$ fixed. In other words, in the flow defined by the initial condition $\phi(x_0,y)=y$. This means that we want to compute
the Taylor expansion coefficients
\[
\varphi_k(x)=\frac {\partial^k \phi}{\partial y^k}(x,y_0)
\]
for which
\[
\phi(x,y)=\varphi(x)+\frac {\partial \phi}{\partial y}(x,y_0)(y-y_0)
+\frac 12\frac {\partial^2 \phi}{\partial y^2}(x,y_0)(y-y_0)^2+\cdots,
\]
is a solution of the equations in variations. The first two variational equations $\VE_1$ and $\VE_2$ are explicitly given by
\begin{equation}\label{VE2}
 \begin{array}{ll}
 \varphi_1^\prime &=f_y(x,\varphi(x))\varphi_1\, , \\
 \varphi_2^\prime &=f_y(x,\varphi(x))\varphi_2+ f_{yy} (x,\varphi(x))\varphi_1^2,
  \end{array}
\end{equation}
where we introduce the standard notation $f_y, f_{yy}, \ldots$ for the corresponding partial derivatives $\partial f /\partial y$, $\partial^2 f/\partial y^2,\ldots$, respectively.
In a similar way we can obtain any higher order variational equation $\VE_k$, for $k>2$.

It is well-known (due to their triangular-shape scheme), that variational equations $\VE_k$ can be solved recurrently. However, in order to apply the Differential Galois Theory, it is convenient to linearise them.
This can be easily done by introducing suitable new variables. For example, if one defines $\chi_1=\varphi_1^2$ and $\chi_2=\varphi_2$, the second variational $\VE_2$ becomes
\begin{equation} \label{LVEE2}
\begin{array}{ll}
 \chi_1^\prime &=2 f_y(x,\varphi(x))\chi_1\, , \\
 \chi_2^\prime &=f_y(x,\varphi(x))\chi_2+f_{yy}(x,\varphi(x))\chi_1.
 \end{array}
\end{equation}
Linearised variational equation of order $k$ will be denoted by $\LVE_k$ (see Appendix~C for more details).
Now, since $\LVE_k$ is a linear system,  standard differential Galois theory can be applied. Let $G_k$ stand for the Galois group of $\LVE_k$ and $G_k^0$ for its identity component.

\begin{definition}
\label{def:meromo:integr}
A complex analytic  field $X$ defined over an analytic complex manifold $M$ of dimension $m$ is \emph{meromorphically integrable} if there exist $X=X_1,X_2,...,X_k$ independent pairwise commuting meromorphic vector fields and $f_1,f_2,...,f_l$ independent meromorphic first integrals of these vector fields, and satisfying that $k+l=m$.
\end{definition}
In the particular planar case, meromorphically integrability can come from two possibilities: 
\begin{itemize}
\item[(a)] either $X_2$ is a symmetry vector field of $X_1$, that is $[X_1,X_2]=0$, or  
\item[(b)] there exists a first integral $H=f_1$ of $X_1$.
\end{itemize}
It is clear that the fibration defined by the level sets of a first integral $H$ of the vector field~(\ref{campo}), $H(x,y)=C$,  gives the solution of equation~\eqref{foliacion}.
Then we have the following result:

\bigskip

\noindent\textbf{Theorem A.}
%\label{MRtheorem}
\textit{
Consider a planar analytic vector field $X$ as in~\eqref{campo} and $\Gamma$ an integral curve.
Then, if $X$ is meromorphically integrable in a neighbourhood of $\Gamma$ it follows that, for any $k\geq 1$, the identity component $G_k^0$ of the Galois group of the linear variational equation $\LVE_k$ is abelian.
}
\bigskip

A very interesting particular case of Theorem~A is when the vector field $X$ is meromorphic at infinity, that is, either is holomorphic or it has a pole at $\infty$. Functions of such type are called meromorphic at the extended complex plane. By compatification of the complex plane, it also known that these functions are the rational functions. This is exactly the case when considering vector fields~\eqref{campo} with $P$ and $Q$ polynomials. 

Therefore, we say that either the field~(\ref{campo}) or the associated foliation~(\ref{foliacion}) is \emph{rationally integrable} analogously as, \emph{mutatis mutandis}, by substituting meromorphic by rational in Definition~\ref{def:meromo:integr}.

\bigskip

\noindent\textbf{Theorem B.}
%\label{MRtheorem}
\textit{
Consider a planar polynomial vector field $X$ as in~(\ref{campo}) and let $\Gamma$ be an algebraic invariant curve.
Assume that the first variational equation $\VE_1$ along $\Gamma$ has an irregular singular point at infinity. Then,
if $X$ is rationally integrable then, for any $k\geq 1$, the Galois group $G_k$ of $\LVE_k$ is abelian.
}
\bigskip

Recall that a linear system of differential equations having a pole at the origin can always be written  in the form $x^h Y'(x) = B(x) Y$, where $B(x)$ is a holomorphic matrix at $x=0$ and $h\in \N$. It is said that $x=0$ is an \emph{irregular singular} point of this equation if $h>1$ (see~\cite{Wasow}, for instance). Equivalently, the solutions of this system around $x=0$ exhibit an exponential growth.
If the first variational equation $\VE_1$ has an irregular singular point, the Galois group $G_k$ of any variational equation is connected and hence equal to its identity component, that is  $G_k=G_k^0$.
Furthermore, if the first variational equation $\LVE_1=\VE_1$ has an irregular (respectively, regular) point then all the following variational equations have the same irregular (respect. regular) point since the exponential behavior depends only on the first variational equation and the solution of this equation appears in the all the following variational equations.

Theorem~A (and~B) are commonly applied to prove non-integrability results: {\it if for some $k$ the Galois group $G_k^0$  ($=G_k$ for Theorem~B) is not abelian then the polynomial field $X$ is not meromorphically/rationally integrable}. Its application is based in the following procedure, similar to the one used for Hamiltonian systems. In the case of Theorem~A:
\begin{itemize}
\item [1)] To obtain an invariant curve $\Gamma$: $y-\varphi(x)=0$, of the field $X$.
\item [2)] To compute $G_k$ and to check if the identity component $G_k^0$ is non abelian (for $k \geq 2$).
\item [3)] If for some $k$, $G_k^0$ is not abelian, then the field $X$ is not meromorphically integrable.
\item [4)] If for all $k$, one has that $G_k^0$ is abelian, then the method cannot decide.
\end{itemize}
Analogously, in the case of Theorem B, the procedure reads:
\begin{itemize}
\item [1')] Same as 1) above.
\item [2')] To compute $G_k$ and to check if it is non abelian (again, for $k\geq 2$).
\item [3')] If for some $k$, we have that $G_k$ is not abelian,  then the field $X$ is not rationally integrable.
\item [4')] If for all $k$, $G_k$ is abelian, then the method cannot decide.
\end{itemize}
It is important to point out that the so--called Risch algorithm (and equation)~\cite{RISCH} provides a systematic way to analyse condition 2') (see Section~\ref{scondicioii} for meore details on it).
The scope of this paper is to use a suitable version of Risch-Kaltofen algorithm (see~\cite{KAL}) to detect the non-integrability of planar polynomial vector fields under the hypotheses of Theorem~B.

\begin{remark}
It is also possible to obtain interesting results by means of Theorem A for {\it regular--singular points}. They are related to the non abelianity of the holonomy of the vector field around the curve $\Gamma$, since for regular--singular points the holonomy of the associated foliation is given, at some order, by the monodromy of the variational equation (see, for instance, \cite{GOM}), and the monodromy group is Zariski dense in Galois group. However, in this regular-singular case it is necessary to check the non commutativity of the identity component of the Galois group.

Like it happens for Hamiltonian systems, the Achilles heel of the above procedure is step 1).
The key point of our approach is that, for generic polynomial vector fields on the plane, the \emph{line of infinity} is always an invariant curve $\Gamma$ of the system (see~\cite{IY}).
In the real section of $\C^2$, this is connected with the Poincar\'e compactification (see, for instance,~\cite[Chapter V]{DLLA} and references therein).
\end{remark}
Along this paper it will be assumed some knowledge on Galois theory of linear differential equations, the so-called Picard--Vessiot theory (we refer reader to~\cite{CH,SV}, and references therein,
for two basic monographs). See also~\cite{morales}.
To our knowledge, this the first time that the  method of the variational equations is applied to study the non--integrability of planar polynomial vector fields.

\bigskip

The paper is structured as follows: Section~\ref{teorema} contains the proof of Theorems~A and~B above. Section~\ref{sstructure} is devoted to the analysis of the Galois group of the variational equations, where conditions (H1) and (Hk) are introduced. These hypotheses will play an essential r\^ole in all the paper. In Section~\ref{scondicioii}, condition (Hk) is directly related to the well-known Risch differential equation. 
Variational equations around the line at infinity are considered in Section~\ref{assum}. We provide general families of examples, perform an explicit \emph{ad-hoc} version of the Risch-Kaltofen algorithm, suitable for our study, in Section~\ref{sec: exa}. A conjecture extending the class of first integrals considered in this paper is presented in Section~\ref{sec:conjecture}.
For the completeness of the work, Appendices on some \textsc{Maple} computations, Galois correspondence and third order variational equations have been included. The main result of this paper (see Theorem~\ref{thm:lstructure:k}), which states

\bigskip
 
\noindent\textbf{Theorem.}
\textit{
Under the assumptions (H1) and (Hk) (for some $k\geq 2$) it follows that the polynomial field $X$ (or its corresponding foliation) is not rationally integrable.
}

\bigskip

\section{Proof of Theorems A and B}
\label{teorema}

We prove Theorem~A. Theorem~B will be a consequence of it in a rational context.
The proof is mainly based on two well-known results:
(i) Morales-Ramis-Sim\'o Theorem and
(ii) Ayoul-Zung Theorem.
Indeed,
\begin{theorem}[\emph{\cite{MRS,morales}}]
\label{MRS}
Assume that a complex analytic Hamiltonian system
is  meromorphically integrable  in a neighbourhood of an integral curve $\Gamma$. Then
the identity components $G_k^0$, $k\geq 1$,  of the Galois groups of the linear variational equations, $\LVE_k$,  along $\Gamma $ are
abelian.
\end{theorem}
Here integrability of a Hamiltonian system means that it is \emph{integrable in the Liouville's sense}: existence of a maximal
number of meromorphic first integrals in involution.

There are several variants of Theorem \ref{MRS}, depending on the nature of the singularities of the first integrals. If the variational equations have irregular singular points at infinity in the phase space then one can compactify the phase space and, therefore, one gets that the obstructions for the existence of (meromorphic) first integrals appear only at the infinity. That is,
\begin{theorem}
\label{MRrational}
Consider a complex analytic Hamiltonian vector field. Assume that the first variational equation $\VE_1$ has an irregular singular point at infinity and that
for some $k$, $G_k^0$ is non abelian. Then the Hamiltonian system is not integrable by means of rational first integrals.
\end{theorem}
%Aqui $G^k=G_0^k???$ Referencia abans tambe

The problem is that, in general, equation~\eqref{foliacion} does not define a Hamiltonian system. Thus, the associated field
\begin{equation}\label{campoY}
Y=\frac{\partial }{\partial x}+f(x,y)\frac{\partial}{\partial y},
\end{equation}
is not Hamiltonian on the phase space $(x,y)$. It is worth noticing that the first integrals of the rational field~\eqref{campoY} are the same as the first integrals of the initial polynomial field~\eqref{campo}, since for any $H(x,y)$ one has that $X(H)=PY(H)$. Thus, although the field $Y$ is not Hamiltonian, it is possible to construct a Hamiltonian system over the cotangent fiber bundle $T^*\mathbb{P}_2$ (called the \emph{cotangent lift of $Y$})\footnote{Since we are assuming that the variational equation has an irregular singular point at infinity, we have to compactify $\C^2$ and, consequently, this must be done in the projective complex plane.} and apply the previous argument on this new Hamiltonian system.
Ayoul and Zung studied this case in~\cite{AZ} and obtained the following more general result:

\begin{theorem}[\emph{\cite{AZ}}] \label{ayoul-zung}
Assume that a complex meromorphic field $Y$ is meromorphically integrable in a neighbourhood
of an integral curve $\Gamma$. Then, for any $k\geq 1$, the identity component $G_k^0$ of the
Galois group of the linear variational equation $\LVE_k$ along $\Gamma $ is abelian.
\end{theorem}
The proof of this theorem is based essentially on the following two facts:
\begin{itemize}
\item The Galois group of the variational equation of the system $Y$ is the same as the Galois group of the variational equations of the lifted system.
\item The application of Theorem~\ref{MRS}, which ensures that integrability of $Y$ implies integrability of the lifted system (in the Liouville sense).
\end{itemize}
Hence, Theorem~\ref{ayoul-zung} undergoes the same variants as Theorem~\ref{MRS}. In particular, for variational equations with irregular singular point at the infinity we have obstructions to the
existence of rational first integrals, i.e. \emph{mutatis mutandis},  Theorem~\ref{MRrational} is also valid here. It is said that the field $Y$ is
{\it rationally integrable} and, therefore, it has a rational first integral.
This concludes the proof of Theorems~A and B.

\begin{remark}
An alternative proof of these theorems can certainly be done using some relevant results by Casale which connect Theorems~A and B with Malgrange's
approach to Galois theory of non-linear differential equations (\cite{CA09}).
\end{remark}
%%%%%%%%%%%%%%%%%%%%%%%%%%%%%%

\section{The Galois group of the variational equation}
\label{sstructure}
Let us focus our attention in the rational context, that is, assume that our vector field is
\[
X=P\dfrac{\partial }{\partial x}+Q \dfrac{\partial }{\partial y},
\]
with $P,Q$ polynomials or, in other words, its foliation is $y'=Q/P=f(x,y)$ with rational $f(x,y)$.
Let us also assume that our invariant curve $\Gamma: y-\varphi(x)=0$ is rational. Then the variational equation \eqref{LVEE2} reads
\begin{equation}
\label{LVE2}  %\tag{$\LVE_2$}
\begin{array}{cl}
\chi_1'&=2\alpha(x)\chi_1\\
\chi_2'&=\alpha(x)\chi_2+\beta(x)\chi_1,
\end{array}
\end{equation}
where $\alpha(x)=f_y(x,\varphi(x))$ and $\beta(x)=f_{yy}(x,\varphi(x))$ are
rational functions. This linear equation~\eqref{LVE2} can be explicitly solved. Indeed, its solutions are given by
%\begin{equation}
\begin{eqnarray}
&&\chi_1(x)= c_1 {\rme}^{2 \int \alpha(x) \, dx}=c_1\omega^2 \\
&&\chi_2(x) = {\rme}^{\int \alpha(x) \, dx} \left( c_1 \int \beta(x) {\rme}^{\int \alpha(x) \, dx} \, dx + c_2 \right)=\omega(c_1\theta_1+c_2), \label{solvar}
\end{eqnarray}
%\end{equation}
with $c_1,c_2 \in \C$ arbitrary constants and provided we define
\begin{equation}
\label{def:omega:theta1}
\omega={\rme}^{\int \alpha}, \qquad \theta_1 = \int \beta \omega.
\end{equation}
Assume that equation~\eqref{foliacion} and the integral curve $\Gamma$ satisfy the following hypotheses:

\medskip

\begin{itemize}
\item [(H1)] Either the rational function $\alpha(x)=R(x)/S(x)$ has a pole (at the finite complex plane) of order greater than 1 or  $\deg(R)\leq \deg(S)$. Moreover, the residues of $\alpha(x)$ at its poles must be all integers.
\item [(H2)] The function $\theta_1$ is not a rational function in the variables $x, \omega$.
\end{itemize}

\medskip

We remark that the first part of hypothesis \textrm{(H1)} is equivalent to the fact that the variational equation~\textrm{$\VE_1$}, $\chi_1'=\alpha(x) \chi_1$, has at least
an irregular singularity in the extended complex plane, that is, a point in which vicinity the solution $\chi_1(x)$ displays an exponential behavior of a rational function. The assumption about the residues is  necessary in order that the integrand  in~\eqref{solvar}, namely,
\[
\beta(x) {\rme}^{\int \alpha(x) \, dx},
\]
could be expressed as a product of a rational function by the exponential of a rational function.

\medskip

From now on we will take, as the field of coefficients of the linear variational equations $\LVE_k$, the rational functions $\C(x)$. Remind that $G_k$ denotes the Galois group of $\LVE_k$, that is $G_k=\textrm{Gal}(\LVE_k)$, and that $G_k^0$ is its corresponding identity component. Thus, the following result holds:
\begin{proposition}
\label{lstructure}
Under assumptions (H1) and (H2) one has that $G_2=G^0_2$ and, moreover, $G_2$ is not abelian.
\end{proposition}
\proof  
A fundamental matrix of system $\LVE_2$ is
\begin{equation}
\label{eq:Phi2}
\Phi_2=\left(
\begin{array}{cc}
\omega^2 &  0\\
\omega\theta_1 & \omega
\end{array}
\right).
\end{equation}
Thus the Galois group $G_2$ is contained in the algebraic group
\[
B:=\left\{ \left(
\begin{array}{cc}
\lambda^2 & 0\\
\lambda\mu  &  \lambda
\end{array}
\right): \lambda\in\C^*, \mu\in\C  \right\},
\]
since, for any $\sigma\in G_2$, we have that
$$\sigma(\omega)=\lambda \omega,\, \, \sigma(\theta_1)=\lambda\theta_1 + \mu, $$
for $\lambda\in\C^*$ (by hypothesis (H1)) and $\mu\in\C $ (by (H2)).
Recall that the dimension of the Galois group is equal to the \emph{transcendence degree} of the Picard--Vessiot extension~\cite{KOL}.
Since the Picard--Vessiot extension corresponding to $\LVE_2$ is
\[
\C(x)\subset \C(x, \omega)=L_1\subset\C\left(x,\omega, \theta_1 \right)=L_2,
\]
one has, from assumption (H1), that
\[
\mathrm{dim}\, G_1=\mathrm{tr\,deg} \left(L_1/\C(x)\right)= \mathrm{tr\,deg}\left(\C(x, \omega)/ \C(x) \right)=1.
\]
So $G_1\simeq \C^*$.
Now assumption (H2) implies that $\theta_1 \not \in \C(x, \omega)$,  so  the transcendence degree $\mathrm{tr\,deg}(L_2/L_1)=1$. This means that the transcendence degree of the Picard--Vessiot extension $\C(x)\subset L_2$ is 2 and the dimension of the Galois group $G_2$ is also 2. Since $G_2$ is contained in $B$, they must coincide, that is $G_2=B$. Moreover, for the algebraic group $B$ we have that $B\simeq \C^* \ltimes \C $ which implies  that $B$ is connected and non abelian. In conclusion, $G^0_2=G_2=B$ is non abelian.

\qed

Hence, assuming that one of the irregular singular points is at the infinity, by  theorem B we obtain the following corollary.
\begin{cor}
\label{cor:lstructure}
Under the assumptions (H1) and (H2) it follows that the polynomial field $X$ (or its corresponding foliation) is not rationally integrable.
\end{cor}
Now, using the Galoisian correspondence, it is possible to generalise Proposition~\ref{lstructure} and Corollary~\ref{cor:lstructure} to higher order variational equations.
Indeed, for $k\geq 2$, we denote by
\[
\beta_k(x)=f_{y\cdots y}(x,\varphi(x)) \, \text{($k$-times derivatives)}, \qquad \theta_{k-1}(x)=\int \beta_k(x)\omega^{k-1}(x)\,dx.
\]
Under this definition the functions $\alpha(x)$, $\beta(x)$ in system~\eqref{LVE2}, read $\alpha(x)=\beta_1(x)$ and $\beta(x)=\beta_2(x)$.
Then, consider the following assumption (which generalises hypothesis (H2)):

\medskip

\begin{itemize}
\item [(Hk)]  The function $\theta_{k-1}$ is not a rational function in the variables $x, \omega$.
\end{itemize}

\medskip

\noindent Therefore, the generalised versions of Proposition~\ref{lstructure} and Corollary~\ref{cor:lstructure} are as follows:
\begin{proposition}
\label{lstructure:k}
Under assumptions (H1) and (Hk) (for some $k\geq 2$) one has that $G_k=G^0_k$ and, moreover, $G_k$ is not abelian.
\end{proposition}

\begin{theorem}
\label{thm:lstructure:k}
Under the assumptions (H1) and (Hk) (for some $k\geq 2$) it follows that the polynomial field $X$ (or its corresponding foliation) is not rationally integrable.
\end{theorem}
The proof of this Theorem follows as a consequence of Proposition~\ref{lstructure:k} and Theorem B.

\medskip

\noindent\textsf{Proof of Prop.~\ref{lstructure:k}.}
Instead of using the explicit linear representation of the Galois group $G_k$ it is more convenient to employ the Galoisian correspondence.
So, the variational equation $\VE_k$ has the following structure
\begin{equation}
\label{VEk}  %\tag{$\LVE_2$}
\begin{array}{rcl}
\varphi_1' &=& \beta_1\varphi_1\\
\varphi_2' &=& \beta_1 \varphi_2 + \beta_2 \varphi_1^2 \\
\varphi_3' &=& \beta_1 \varphi_3 + 3 \beta_2 \varphi_1 \varphi_2 + \beta_3 \varphi_1^3 \\
\vdots     &\vdots& \vdots \\
\varphi_k'&=& \beta_1 \varphi_k+ \cdots + \beta_k(x)\varphi_1^k
\end{array}
\end{equation}
and its corresponding $\LVE_k$ becomes
\begin{equation}
\label{LVEk}  %\tag{$\LVE_2$}
\begin{array}{rcl}
\chi_1'&=&k\beta_1 \chi_1\\
\vdots &\vdots& \vdots \\
\chi_k'&=& \beta_1 \chi_k + \cdots + \beta_k(x)\chi_1,
\end{array}
\end{equation}
with $\chi_1=\varphi_1^k, \ldots,\chi_k=\varphi_k$ and where only relevant terms have been explicitly written.
From this structure and as it happens for $\LVE_2$, we see that the Picard--Vessiot extension $\C(x)\subset L_k$ of $\LVE_k$ is obtained by means of $\omega$ (exponential of a primitive) following
some quadratures. Hence, like for $k=2$, assumption (H1) implies that the Galois group $G_k$ must be connected.

Since $\chi_1=\omega^k$, and $\chi_k=c_1\omega+c_2\omega\theta_{k-1}+\cdots $, the solution of $\LVE_k$ must contain $\omega$ and  $\theta_{k-1}$. Therefore, one has the intermediate extension
\[
\C(x)\subset \C(x,\omega,\theta_{k-1})\subset L_k.
\]
We note that  the extension $K=\C(x)\subset \C(x,\omega,\theta_{k-1})=S$ is a  Picard--Vessiot extension of the linear subsystem
\begin{equation}
%y^\prime+(k-1)\alpha y=\beta_k.
\begin{array}{rcl}
\bar{\chi}_1'&=&k\beta_1 \bar{\chi}_1\\
\bar{\chi}_k'&=& \beta_1 \bar{\chi}_k + \beta_k(x)\bar{\chi}_1,
\end{array}
\label{rischk}
\end{equation}
Let $H=Gal(L_k/S)$ the Galois group of the extension $S\subset L_k$ and $G=Gal(S/K)$ the Galois group of equation~\eqref{rischk}. From hypotheses (H1) and (Hk) it is easy to prove that the group $G$ is not abelian (in fact, as for the second order variational equation,  it is a semidirect product of the additive by the multiplicative group). So, by the Galoisian correspondence of the Picard--Vessiot theory (we refer the reader to Appendix~B for more details), $G\simeq G_k/H$.
As $G$ is not abelian it follows that $G_k$ is not abelian as well.

\qed

\begin{remark}
In Appendix~C the representation of the Galois group of the third variational equation is explicitly given. The reader can easily derive the proof of Proposition~\ref{lstructure:k} for the case $k=3$.
\end{remark}

The results above draw a systematic scheme to check the non-rational integrability of a given polynomial field: first verifying (H1); then ckecking (H2) and if it is not satisfied go for (H3) and so on. If we find a $k\geq 2$ such that (Hk) holds, then our field is not rationally integrable.

As it will be seen later, equation~\eqref{rischk} (which appears in a natural way in the previous proof) will play a crucial r\^ole in the rest of the paper. In fact, when facing a concrete family of polynomial systems, to check hypothesis (H1) will be quite straightforward by direct inspection of the rational function $\alpha(x)$. However to check whether hypothesis (H2) (or more general (Hk)) is satisfied is a much more involved difficulty. Next section is devoted to the relation between this problem and the existence of a rational solution of the so-called Risch equation.

\section{Liouville's Theorem and Risch differential equation}
\label{scondicioii}

Let $f(x)$ and $g(x)$ be rational functions with $g(x)$ non constant.
We say that the integral
${\ds \int f(x){\rme}^{g(x)}\,dx}$ is \emph{elementary} if it can be expressed in the form
\[
\int f(x){\rme}^{g(x)}\,dx = {\rme}^{g(x)} h(x)+c,
\]
where $h(x)$ is a rational function and $c$ is an integration constant.
The following assertion is a particular result from a more general theorem of Liouville.
\begin{theorem}[\cite{LIO}]
\label{strong}
%\label{strong:lio}
The integral ${\ds \int f(x){\rme}^{g(x)}\,dx}$ is elementary if and only if the differential equation
\begin{equation}
\label{RISCH}
y'+ g'(x) y=f(x)
\end{equation}
has a rational solution $h(x)$.
More precisely, one has that ${\ds \int f(x){\rme}^{g(x)}\,dx = {\rme}^{g(x)} h(x) }+ \textrm{constant}$ if and only if $h(x)$ satisfies  equation \eqref{RISCH}.
\end{theorem}
Equation~\eqref{RISCH} is usually called the \emph{Risch differential equation}, since he was the first one who provided an algorithm to decide whether it has or not a rational solution (see~\cite{RISCH}). For a proof of this theorem using differential fields we refer the reader to~\cite[p.46]{RITT}.

In our setting we have $f(x)=\beta(x)$ and $g(x)=\int \alpha(x) \, dx$. So, Risch equation~\eqref{RISCH} becomes
\begin{equation}
\label{RISCHii}
y'+\alpha(x)y=\beta(x).
\end{equation}
Hypothesis (H1) implies that $\omega$ can be reduced to a transcendental function
of the form $h_1(x){\rme}^{g_1(x)}$, with $h_1$ and $g_1$  rational functions. The function $h_1$
comes from the cancelation of the exponential with the possible logarithms defined by the  poles of $\alpha$ with integer residues. Thus,
in the case that the function $\alpha$ has  terms in $1/(x-x_i)$, that is,
 $${\ds \alpha(x) =\tilde{\alpha}(x) + \sum_i \frac{\ell_i}{x-x_i}}, $$
  with $\ell_i\in \Z$, then equation\eqref{RISCHii} must be changed by
\[
y'+\tilde{\alpha}(x)y=\tilde{\beta}(x),
\]
where $\tilde{\beta}(x)=\underset{i}\prod(x-x_i)^{\ell_i} \beta(x)$.

On the other side, hypothesis (H2) is equivalent to say that the Risch equation~\eqref{RISCHii} admits no rational solution. It is worth to mention that if (H1) applies and (H2) fails then there exists a unique rational solution $h(x)$ of the Risch equation. Indeed,
the general solution of equation~\eqref{RISCH} is  $y(x)=c\, {\rme}^{-g(x)} + {\rme}^{-g(x)} \int f(x) {\rme}^{g(x)}\, dx$, $c\in \C$ a general constant. Since (H2) fails then $\int f(x){\rme}^{g(x)}\,dx$ is elementary and so $\int f(x){\rme}^{g(x)}\,dx = {\rme}^{g(x)} h(x)$ and the general solution becomes $y(x)=c{\rme}^{-\int \alpha} + h(x)$, which is rational only for $c=0$.

Hence, by Theorem~\ref{strong}, assuming that hypothesis (H1) is satisfied, condition (H2) to be fulfilled is equivalent to the fact that the
Risch equation~\eqref{RISCHii} has no rational solution. The extension to (Hk) is straightforward since its associated Risch equation is given by
\begin{equation}
\label{eq:Risch:order:k}
y'+(k-1)\alpha(x) y = \beta_k(x). 
\end{equation}
For algorithms dealing with this problem we refer the reader to the references~\cite{DAV,ROT,KAL}.
We stress that there are usually two ways to approach this kind of issues:
\begin{itemize}
\item [a)]  \emph{Analytic approach}, which tries to prove by analytic methods whether the equation~(\ref{RISCHii}) (or more general \eqref{eq:Risch:order:k}) has or not any rational solution.

\item [b)]  \emph{Algebraic approach}, seeking to prove or disprove in a direct way the existence (or not) of a rational solution to equation \eqref{RISCHii} (or more general, for~\eqref{eq:Risch:order:k}).
This is the way employed, essentially, in this work.

\end{itemize}

\section{Variational equations along the line at infinity}
\label{assum}

In this section we consider polynomial systems and study their variational equations along the line at infinity.
We follow the ideas introduced in~\cite{GOM}.
Let consider the polynomial field:
\begin{equation}
\widetilde{X}=\widetilde{P}(z_1,z_2) \dfrac{\partial }{\partial z_1}+\widetilde{Q}(z_1,z_2) \dfrac{\partial }{\partial z_2},
\label{Vectorz}
\end{equation}
with
\begin{equation}
\label{prop:ve:polPQ}
\widetilde{P}(z_1,z_2)=\sum_{i=0}^{N_1} P_i(z_1,z_2), \qquad
\widetilde{Q}(z_1,z_2)=\sum_{i=0}^{N_2}  Q_i(z_1,z_2),
\end{equation}
$P_i(z_1,z_2)$, $Q_i(z_1,z_2)$ being homogeneous polynomials of degree $i$. Let define $N:=\max \left\{ N_1, N_2 \right\}$, the degree of $X$.
The foliation defined by the field $\widetilde{X}$ is given by the first order differential equation
\begin{equation}
\frac {dz_2}{dz_1}=\frac{\widetilde{Q}(z_1,z_2)}{\widetilde{P}(z_1,z_2)}.
\label{fol}
\end{equation}
In case that $\widetilde{P}\equiv 0$ one should interchange the r\^oles of $z_1$, $z_2$ and $\widetilde{P}$,$\widetilde{Q}$.
Two fields are said to be {\it equivalent} if they define the same foliation. Thus, in a dynamical language, we are more concerned with their orbits than with the time parametrisation of the integral curves of the field.
\begin{proposition}
\label{prop:line:y0}
For $z_1\neq 0$, the birational change of coordinates $y=1/z_1$ and $x=yz_2=z_2/z_1$ leads
the polynomial vector field $\widetilde{X}$ to the rational system
\[
\begin{array}{rcl}
\dot{x}&=& \phantom{+}\displaystyle{\sum_{i=0}^Ny^{1-i}\left( Q_i(1,x)-xP_i(1,x) \right)}\\
\dot{y}&=& -\displaystyle{\sum_{i=0}^Ny^{2-i}P_i(1,x)}.
\end{array}
\]
This change of variables sends the line at infinity $z_1=\infty$ to the line $y=0$. Its corresponding foliation is defined by
the rational ODE
\begin{equation}
\label{tildef}
\dfrac{dy}{dx}=\dfrac{\displaystyle{y\sum_{i=0}^Ny^{N-i}P_i(1,x)}}{
\displaystyle{\sum_{i=0}^Ny^{N-i} \left( xP_i(1,x)-Q_i(1,x) \right)}}
\end{equation}
and the associated polynomial field is given by
\begin{equation}
\label{prop:line:X}
X= P(x,y) \frac{\partial}{\partial x} + Q(x,y) \frac{\partial}{\partial y},
\end{equation}
where
\begin{equation}
\label{prop:line:PQ}
P(x,y):= \sum_{i=0}^N y^{N-i} \left( x P_i(1,x)-Q_i(1,x)\right), \qquad Q(x,y):=y \sum_{i=0}^N y^{N-i} P_i(1,x).
\end{equation}
Moreover, the following relations between $\widetilde{P}, \widetilde{Q}$ and $P,Q$ hold:
\[
P(x,y) = \frac{1}{z_1^{N+1}} \widetilde{P}(z_1,z_2), \qquad
Q(x,y)= \frac{1}{z_1^{N+1} z_2} \left( \widetilde{P}(z_1,z_2) + z_1 \widetilde{Q}(z_1,z_2) \right),
\]
and
\[
\widetilde{P}(z_1,z_2) = \frac{1}{y^{N+1}} P(x,y), \qquad
\widetilde{Q}(z_1,z_2) = \frac{1}{y^{N+1}} \left( -y P(x,y) + x Q(x,y) \right).
\]
\end{proposition}

\proof
From the change of variables it is clear that
\begin{eqnarray*}
\dot{x} &=& \frac{1}{z_1} \dot{z}_2 - \frac{z_2}{z_1} \, \frac{\dot{z}_1}{z_1} =
\frac{1}{z_1} \left( \sum_{i=0}^N z_1^i Q_i(1,z_2/z_1) \right) - \frac{z_2}{z_1^2} \left( \sum_{i=0}^N z_1^i P_i(1,z_2/z_1) \right) = \\
&& \sum_{i=0}^N y^{1-i} \left( Q_i(1,x) - x P_i(1,x) \right).
\end{eqnarray*}
Analogously,
\[
\dot{y} = - \frac{\dot{z}_1}{z_1^2} = - y^2 \sum_{i=0}^N y^{-i} P_i(1,x) = - \sum_{i=0}^N y^{2-i} P_i(1,x).
\]
Notice that, however the polynomials $P(z_1,z_2)$ and $Q(z_1,z_2)$ were homogeneous in $z_1,z_2$
$P_i(1,x)$ and $Q_i(1,x)$, polinomials of degree at
most $i$, are in general nonhomogeneous.
The corresponding foliation of the system above is given by
\[
\frac{dy}{dx}= \frac{\dot{y}}{\dot{x}} = \frac{y {\ds \sum_{i=0}^N y^{N-i} P_i(1,x)}}{{\ds \sum_{i=0}^N y^{N-i} \left( x P_i(1,x) - Q_i(1,x) \right)}} =:
\frac{Q(x,y)}{P(x,y)},
\]
where we have multiplied $\dot{y}$ and $\dot{x}$ by $y^{N-1}$. Concerning the relations between $\widetilde{P},\widetilde{Q}$ and $P,Q$, one has
\[
\widetilde{P} (z_1,z_2) = \sum_{i=0}^N z_1 ^i P_i(1,z_2/z_1) = \frac{1}{y^{N+1}} P(x,y)
\]
or, equivalently,
\[
P(x,y)=\frac{1}{z_1^{N+1}} \widetilde{P}(z_1,z_2).
\]
In a similar way one gets
\[
\widetilde{Q}(z_1,z_2) = \sum_{i=0}^N z_1^i Q_i(1,z_2/z_1) = \frac{1}{y^{N+1}} \left( -y P(x,y) + x Q(x,y) \right)
\]
and, taking into account the relation between $\widetilde{P}$ and $P$,
\[
Q(x,y) = \frac{1}{z_1^{N+1} z_2} \left( \widetilde{P}(z_1,z_2) + z_1 \widetilde{Q}(z_1,z_2) \right).
\]
So, in the end,
\begin{equation}
\frac{dz_2}{dz_1}= \frac{\widetilde{Q}(z_1,z_2)}{\widetilde{P}(z_1,z_2)} =
\frac{z_2{Q}(z_2/z_1,1/z_1)-{P}(z_2/z_1,1/z_1)}{z_1{Q}(z_2/z_1,1/z_1)}.
\label{foliz}
\end{equation}

\qed

\medskip
\noindent From now on, we restrict ourselves to polynomial vector fields $X$ of the form~\eqref{prop:line:X} with $P(x,y)$, $Q(x,y)$ as in~\eqref{prop:line:PQ} and having $y=0$ as an invariant curve.
It is straightforward to check in this case that the corresponding second order linear variational equation $\LVE_2$ becomes
\begin{equation}
\label{LVE2:gen}
\begin{array}{rcl}
\chi_1' &=& 2 \alpha(x) \chi_1 \\[1.2ex]
\chi_2' &=& \beta(x) \chi_1 + \alpha(x) \chi_2,
\end{array}
\end{equation}
with
\begin{equation}
\label{LVE2:polynomial}
\begin{array}{lll}
\alpha(x) &=&\dfrac{P_N(1,x)}{xP_N(1,x)-Q_N(1,x)}, \\[1.2ex]\\&&\\
\beta(x)  & =&2\, \dfrac{P_N(1,x)Q_{N-1}(1,x)-P_{N-1}(1,x)Q_N(1,x)}{\left( xP_N(1,x)-Q_N(1,x)\right)^2}.
\end{array}
\end{equation}

\section{Examples and algorithmic considerations}
\label{sec: exa}

To illustrate the use of the method, we will focus our attention on those families of polynomial vector fields $X$ of the form
\begin{equation}
\label{subfamily}
\dfrac{dy}{dx}=\dfrac{y(P_1(x) y^{N-1}+\cdots+P_{N-1}(x)y+P_N(x))}{x^k-y},
\end{equation}
or those families $\widetilde{X}$, as in~\eqref{Vectorz}, that can be led into this form by means of the birational change of variables $y=1/z_1$, $x=z_2/z_1$. To avoid a cumbersome notation, hereafter we will also denote $P_j(x)=P_j(1,x)$.
\begin{proposition}
Let consider integers $N\geq 2$, $2\leq k \leq N$ and a polynomial $\widetilde{P}$ of degree $N$,
\[
\widetilde{P}(z_1,z_2)=\sum_{i=0}^N P_i(z_1,z_2),
\]
with $P_i(z_1,z_2)$ homogeneous polynomial of degree $i$, satisfying that $P_0=0$ and $P_i(0,z_2)=0$, for $i=1,\ldots,N$.
Let us define a polynomial of degree $N$,
\[
\widetilde{Q}(z_1,z_2)=\sum_{i=0}^N Q_i(z_1,z_2),
\]
where $Q_i(z_1,z_2)$, also homogeneous polynomials of degree $i$, are given by the following relations:
\begin{eqnarray}
Q_0(z_1,z_2)      &\equiv& 0, \nonumber \\[1.2ex]
Q_{\ell}(z_1,z_2) &=& {\ds \frac{z_2}{z_1}} P_{\ell}(z_1,z_2), \qquad \ell=1,\ldots, N-2 \nonumber \\[1.2ex]
Q_{N-1}(z_1,z_2)  &=& {\ds \frac{z_2}{z_1}} P_{N-1}(z_1,z_2) + z_1^{N-1} \label{prop:Qfamily} \\[1.2ex]
Q_N(z_1,z_2)      &=& {\ds \frac{z_2}{z_1}} P_N(z_1,z_2) - z_1^{N-k} z_2^k. \nonumber
\end{eqnarray}
Under the transformation $y=1/z_1$, $x=z_2/z_1$, the new polynomial vector field $X$ becomes
\begin{equation}
\label{eq:family:PQ}
X = \left( x^k - y \right) \frac{\partial}{\partial x} + y \left( P_1(x) y^{N-1} + \cdots + P_{N-1}(x) y + P_N(x) \right) \frac{\partial }{\partial y}
\end{equation}
and its foliation takes the form~\eqref{subfamily}.
\end{proposition}

\proof From Proposition~\ref{prop:line:y0} we know that under the transformation $y=1/z_1$, $x=z_2/z_1$, any polynomial vector field $\widetilde{X}$ given by~\eqref{Vectorz} and~\eqref{prop:ve:polPQ} takes the form
\[
\frac{dy}{dx}= \frac{y {\ds \sum_{i=0}^N y^{N-i} P_i(1,x)}}{{\ds \sum_{i=0}^N y^{N-i} \left( x P_i(1,x) - Q_i(1,x) \right)}}.
\]
So, in our case we should determine $\widetilde{Q}(z_1,z_2)$ such that
\[
\sum_{i=0}^N y^{N-i} \left( Q_i(1,x) - x P_i(1,x) \right) = x^k -y.
\]
Equating powers in $y$, one gets $Q_0=P_0=0$ and
\begin{eqnarray*}
Q_{\ell}(1,x) &=& x P_{\ell}(1,x), \quad \ell=1,\ldots N-2, \\
Q_{N-1}(1,x)  &=& xP_{N-1}(1,x)+1, \qquad Q_N(1,x)=xP_N(1,x) - x^k.
\end{eqnarray*}
Expressed in $(z_1,z_2)$-variables, they are led into the expression~\eqref{prop:Qfamily}. Notice that the terms $(z_2/z_1)P_j(z_1,z_2)$ are well-defined polynomials since, by hypothesis, $P_j(0,z_2)=0$, for $j=1,\ldots,N$.

\qed

From now on, we restrict ourselves to polynomial vector fields $X$ of degree $N$ whose foliation is given by a differential ODE of the form~\eqref{subfamily},
\[
\dfrac{dy}{dx}=\dfrac{y(P_1(x) y^{N-1}+\cdots+P_{N-1}(x) y+P_N(x))}{x^k-y}:=f(x,y)
\]
and with $2\leq k \leq N$, an integer. As it was already pointed out at the end of Section~\ref{assum}, its corresponding second order linear variational equation is given by~\eqref{LVE2:gen}. In our case the invariant curve is $\Gamma: \ y=0$, that is $y=\varphi(x)=0$, so
\begin{eqnarray*}
\alpha(x) &=& f_y(x,\varphi(x))= f_y(x,0) = \frac{P_N(1,x)}{x^k} \\
\beta(x)  &=& f_{yy}(x,\varphi(x))=f_{yy}(x,0) = 2 \left( \frac{P_N(1,x)}{x^{2k}} + \frac{P_{N-1}(1,x)}{x^k} \right),
\end{eqnarray*}
as it was already stated in~\eqref{LVE2:polynomial}.

\medskip

\noindent\textbf{Example 1.}
%\begin{example}
%\label{ex:errorf}
Consider the planar vector field
\begin{equation}
\label{errorvf}
X=(x^3-y)\dfrac{\partial}{\partial x}+y(x^2-cx-b-ay)\dfrac{\partial}{\partial y}
\end{equation}
and its associated foliation defined by
\begin{equation}\label{errorf}
\dfrac{dy}{dx}=\dfrac{y(x^2-cx-b-ay)}{x^3-y}.
\end{equation}
Note that the field $X$ has invariant the straight line $y=0$. Now we consider the variational equations
$\LVE_2$ of \eqref{errorf} on the line $y=0$ and we obtain
\begin{equation}
\begin{array}{cl}
\chi_1^\prime=&2\ \dfrac{x^2-cx-b}{x^3}\chi_1\\
\chi_2^\prime=&\dfrac{x^2-cx-b}{x^3}\chi_2-2\left(\dfrac{a}{x^3}-\dfrac{x^2-cx-b}{x^6}
\right)\chi_1.
\end{array}
\end{equation}
Condition $b\neq 0$ implies hypothesis (H1). In order to check hypothesis (H2) we reduce our problem to study if its associated Risch equation
\begin{equation} \label{eqdif:fam}
y^\prime+\dfrac{x^2-cx-b}{x^3}y=-
2\left(\dfrac{a}{x^3}-\dfrac{x^2-cx-b}{x^6}\right)
\end{equation}
has or not a rational solution.
%\end{example}

\bigskip

\noindent\textbf{Example 2 (an infinite family).} 
%\label{ex:infin} 
Assume that in~\eqref{subfamily} $P_N=a\in\mathbb{C}$ and $P_{N-1}=b\in\mathbb{C}$.
We assume $a\neq 0$ and $k\geq 2$. The field is
\begin{equation}
\label{eq:infinity}
X=(x^k-y)\dfrac{\partial}{\partial x} +
y(P_1(x)y^{N-1}+\cdots P_{N-2}(x)y^2+by+a)\dfrac{\partial}{\partial y}
\end{equation}
with foliation
\begin{equation}
\label{ex2:errorf}
\dfrac{dy}{dx}=\dfrac{y(P_1(x)y^{N-1}+\cdots P_{N-2}(x)y^2+by+a)}{x^k-y},
\end{equation}
$b\in \C, a\in\C^*.$ Along the straight line $y=0$ its $\LVE_2$ becomes
\[
\begin{array}{cl}
\chi_1'&={\ds 2\frac{a}{x^k}\chi_1}\\[1.2ex]
\chi_2'&={\ds \frac{a}{x^k}\chi_2+2\left(\frac{a}{x^{2k}}+\frac{b}{x^k}\right)\chi_1}.
\end{array}
\]
As before, hypothesis (H2) and, therefore its rational integrability, can be reduced to study the existence of a rational solution of the corresponding Risch equation
\begin{equation} 
\label{ex2:eqdif:fam}
y^\prime+\frac a{x^k}y=\frac {2a}{x^{2k}}+\frac {2b}{x^k}.
\end{equation}
We will prove later that this is not the case for $k>2$.
%\end{example}

\bigskip

\noindent\textbf{Example 3.}
%\begin{example} 
%\label{ex:whit}
Let $X$ be the following quadratic polynomial field:

\begin{equation}\label{whit}
X=y(a_1x+a_0)\dfrac{\partial}{\partial y}+(x^2-y)\dfrac{\partial}{\partial x}.
\end{equation}
Then
\begin{equation}
\label{whitf}
\dfrac{dy}{dx}=\dfrac{y(a_1x+a_0)}{x^2-y}.
\end{equation}
Despite its simple appearance we will prove that generically the above field is not rationally integrable.
Along the straight line $y=0$ the LVE$_2$ becomes
\[
\begin{array}{cl}
\chi_1'&={\ds 2\frac{a_1x+a_0}{x^2}\chi_1} \\[1.2ex]
\chi_2'&={\ds \frac{a_1x+a_0}{x^2} \chi_2+ \frac{2a_1x+2a_0}{x^4}\chi_1}.
\end{array}
\]
As before, we study the existence of rational solutions for its corresponding Risch equation:
\[
y'+ \frac{a_1 x+a_0}{x^2} y = 2\, \frac{a_1 x+a_0}{x^4}.
\]

\bigskip

\subsection{Application of the Risch-Kaltofen algorithm}
\label{sec:alg}

In this section the solution of the Risch equation~\eqref{RISCHii} is analysed by means of the so--called Risch algorithm. It follows the ideas of the work by Kaltofen in~\cite{KAL} (from now on refered as Risch-Kaltofen algorithm). We restrict ourselves here to the particular case where
$\alpha$ and $\beta$ are given by the expressions
\begin{equation}
\label{alg}
\begin{array}{lcl}
{\ds \alpha(x)=\frac{A(x)}{x^{k}} }, &\quad & {\ds \beta(x)={ \frac{2A(x)+2x^kB(x)}{x^{2k}}} }, \\[1.2ex]
{\ds A(x)=\sum_{i=0}^{n}a_{i}x^i}, &\quad& {\ds B(x)=\sum_{i=0}^{m}b_{i}x^i,\quad k\in \mathbb{Z}^+},
\end{array}
\end{equation}
with $a_n\neq 0$, $b_m\neq 0$, $k>1$ and $n<k$. From now on we will always consider $a_0\neq 0$ (if not, this would imply that $k$ becomes, at least, $k-1$ in $\alpha(x)$).  
It is clear that  $\alpha(x)$ has a pole of order greater than 1 at $x_0=0$ (hypothesis (H1)).

\begin{remark}
Relevant special cases of the family~\eqref{eq:family:PQ} can be expressed in the latter form. For instance,
\[
\alpha(x)=\frac{P_N(x)}{x^k}, \qquad \beta(x)=2 \left( \frac{P_N(x)}{x^{2k}} + \frac{P_{N-1}(x)}{x^k} \right),
\]
with $n=\deg P_N(x)$, $m=\deg P_{N-1}(x)$, $n<k$. In particular,  for $k=N+1$, we have $n=\deg P_N(x)\leq N<k$.  However, as we will see later,   there are polynomial systems not included in the family~\eqref{eq:family:PQ} with $\alpha$ and $\beta$ as in~\eqref{alg}.
\end{remark}
Observe that $a_0 \ne 0$ implies that $\deg(\mathrm{GCD}(A(x),x^k))=0$ and also that $\deg(\mathrm{GCD}(2A(x)+2^kB(x),x^{2k}))=0$ (since they are irreducible fractions). In the same way, $b_0\neq 0$, unless $B(x)$ be identically zero, i.e., $B\equiv 0$. 
According to the notation of Kaltofen \cite{KAL}, 
\[
q_1=x,\, k_1=k,  \, l_1=2k,\, F=A(x),\, G=2A(x)+2x^kB(x).
\] 
Due to $k>1$, then 
\[
\tilde{j}_1=\min\{2k+1,k\}=k,\quad y(x)=\frac{Y(x)}{x^k} 
\]
and we arrive to the following differential equation $uY'+vY=w$, i.e.,
\begin{equation}\label{eqka}
x^kY'+(A-kx^{k-1})Y=2A+2x^kB,\end{equation} where   \begin{equation}\label{eqka2}\begin{array}{ll}
u=u_px^p+\ldots+u_0=x^k,& v=va_rx^r+\ldots+v_0=A-kx^{k-1},\\
w=w_sx^s+\ldots+w_0=2A+2x^kB,& Y=y_hx^h+\ldots+y_0.
\end{array}
\end{equation} 
Replacing $x=0$ in equation \eqref{eqka} and having in mind that $a_0\neq 0$, it follows that
$y_0=2$. Since  $n<k$, $s<2k$ it follows $m<k$ and $n\leq k-1$, and thus we arrive to
\[
p=k,\quad u_p=1,\quad r=k-1, \quad s=\left\{\begin{array}{ll} m+k,& m\geq 0\\ n,& B\equiv 0.
\end{array}\right.
\]
Let us define
\[
v_r=\left\{\begin{array}{ll}
-k,& n<k-1\\
a_{n}-k, & n=k-1,
\end{array}\right.
\]
and the function $\rho=-v_r/u_p$ as
\[
\rho=\left\{\begin{array}{cl}
-v_r& v_r\in \mathbb{Z}^-\\
0 & v_r\notin \mathbb{Z}^-.
\end{array}\right.
\]
Therefore, equation \eqref{eqka} can be rewritten as
\begin{equation}\label{eqka3}
\sum_{i=1}^{\widetilde{h}}\left[(i-k)x^{i+k-1}+Ay_ix^i\right]=\sum_{i=0}^{m}\left(2b_ix^{k+i}\right)+2kx^{k-1},
\end{equation}
where $\widetilde{h}\geq h$ is given by
\begin{equation*}
%\label{eqwih}
\widetilde{h}=\left\{\begin{array}{ll}
\max\{\min\{m,m+1\},\rho\}                         & \textrm{if $m\geq 0$,}\\
\max\{\min\{n-k-1,n-k+1\},\rho\}=\max\{n-k-1,\rho\}& \textrm{if $B\equiv 0$.}
\end{array}
\right.
\end{equation*}
In our case, this $\widetilde{h}$ reads
\begin{equation}
\label{eqwih}
\widetilde{h}=\left\{\begin{array}{ll}
\max\{m,\rho\} & \textrm{if $m\geq 0$,}\\
\rho & \textrm{if $B\equiv 0$.}
\end{array}
\right.
\end{equation}
One should obtain the rest of coefficients of $Y(x)$, that is $y_{\widetilde{h}}, y_{\widetilde{h}-1},\ldots,y_1$, whenever the algebraic equation \eqref{eqka} has solution. Applying Rouch\'e-Frobenius theorem over the system obtained after specialisation of the curve with $\widetilde{h}$ points $x_1,\ldots x_{\widetilde{h}}$ (where $x_i\neq 0$ for all $1\leq i\leq \widetilde{h}$) the polynomial $Y(x)$ exists if and only if the rank of the matrix of the system is exactly the rank of the augmented matrix. Recall that if
\begin{equation}\label{eqka4}
\deg(x^kY'+(A-kx^{k-1})Y)\neq\deg(2A+2x^kB),
\end{equation}
then there is no solution for $Y(x)$.

According to possibilities for $\widetilde{h}$ in equation \eqref{eqwih} we have the following two cases:
\begin{enumerate}
\item Case 1: $n<k-1$. This is the trivial case because $v_r=-k\in\mathbb{Z}^-$ and since to $m<k$, it derives, by equation~\eqref{eqwih}, that $\widetilde{h}=k$. By~\eqref{eqka3} we have $$\sum_{i=1}^{k}\left[(i-k)x^{i+k-1}+Ax^i\right]y_i=\sum_{i=0}^{m}\left(2b_ix^{k+i}\right)+2kx^{k-1},$$
which, after cancellation of the term of degree $2k-1$, lead us to 
\[
\left\{\begin{array}{ll}
-y_{k-1}x^{2k-2}+\ldots+(1-k)y_1x^k+a_ny_kx^{n+k}+\ldots+ a_0y_1x\equiv 0,& B\equiv 0,\\
-y_{k-1}x^{2k-2}+\ldots+(1-k)y_1x^k+a_ny_kx^{n+k}+\ldots+ a_0y_1x=2b_mx^{m+k}+\ldots 2b_0x^k,& m\geq 0
\end{array}\right.
\]
Observe that for $B \equiv 0$, all the coefficients on the left-hand side must vanish in order to provide solution for $Y(x)$. On the other hand, for $m\geq 0$, the degree of the polynomial in the left-hand side is at most $2k-2$ and the degree of the polynomial in the right-hand side is at most $m+k<2k$. So, we can obtain rational solutions of Risch differential equation when equation \eqref{eqka3} has solution for $Y(x)$. In particular, by condition~\eqref{eqka4}, if the degree of the left-hand side is $2k-2$, the degree of the right-hand side is $m+k$ and since $m\neq k-2$ the Risch differential equation has no rational solution.

\item Case 2: $n=k-1$
\begin{enumerate}
    \item If $m<1$ and $k-a_n\notin \mathbb{Z}^+$, then $\rho=0$ and by equation~\eqref{eqwih} we obtain $\widetilde{h}=0$. From~\eqref{eqka2} we get $Y(x)=y_0=2$. Replacing $Y(x)=2$ into equation~\eqref{eqka3} we obtain $xB=-k$, which is a contradiction because it does not exist $m\in\mathbb{N}$ such that $xB=-k$ (although $B\equiv 0$, because $k\neq 0$). Thus, we conclude that in this case there are no rational solutions for the Risch differential equation.
    \item If $k-a_n\notin \mathbb{Z}^+$ and $m\geq 1$, then $\rho=0$ and by equation \eqref{eqwih} we obtain $\widetilde{h}=m$. Using~\eqref{eqka2} we see that $Y(x)=y_{m}x^{m}+\ldots +2$ and replacing it into equation~\eqref{eqka3} we get rational solutions of Risch differential equation provided the rest of indeterminate coefficients $y_1,\ldots,y_{m}$ satisfy the algebraic equation \eqref{eqka}.
    \item If $k-a_n\in \mathbb{Z}^+$ and $m\geq k-a_n$, then $\rho=a_n-k$ and by equation~\eqref{eqwih} we obtain $\widetilde{h}=m$. Similarly, from~\eqref{eqka2} we see that $Y(x)=y_{m}x^{m}+\ldots +2$, and replacing it into equation~\eqref{eqka3}, we obtain rational solutions of the Risch differential equation provided the rest of indeterminate coefficients $y_1,\ldots,y_{m}$ satisfy the algebraic equation~\eqref{eqka}.
    \item If $k-a_n\in \mathbb{Z}^+$ and $m<k-a_n$, then $\rho=a_n-k$ and by equation \eqref{eqwih} we obtain $\widetilde{h}=k-a_n$. Taking into account~\eqref{eqka2} one sees that $Y(x)=y_{k-a_n}x^{k-a_n}+\ldots +2$ and by replacing it into equation~\eqref{eqka3} one gets rational solutions for the Risch differential equation whenever the rest of indeterminate coefficients $y_1,\ldots,y_{k-a_n}$ satisfy the algebraic equation~\eqref{eqka}.
\end{enumerate}
    \end{enumerate}

\medskip    
    
\noindent Now we go back to our families of examples.

\medskip

%\begin{example}
%\label{ex:errorf}
\noindent\textbf{Example 1 (continuation).}
Let us take $k=3$, $n=2$, $m=0$, $a_2=1$, $a_1=-c$, $a_0=-b$, $b_0=-a$ in equation~\eqref{errorvf}. Therefore:
\begin{itemize}
\item In the case $a_1 \neq a_0b_0/3$ (equivalently, $c\neq -ab/3$), it follows that $n=k-1$ and $m=0$ so we fall in case 2. Precisely, since $v_r=a_2-k=-2\in \mathbb{Z}^-$, we fall in subcase 2.d, which lead us to $\widetilde{h}=2$. Hence, $Y(x)=y_2x^2+y_1x+2$. Now, substituting $Y(x)$ into equation \eqref{eqka3} we see that condition \eqref{eqka4} holds and Risch differential equation has no rational solution. In consequence, equation~\eqref{errorvf} is no rationally integrable.
Coming back to the field $\widetilde{X}(z_1,z_2)$, by means of equation~\eqref{foliz}, we obtain that the field
\begin {equation}
\widetilde{X}=(ab z_1^2+cz_1z_2-z_2^2+z_1)\frac{\partial}{\partial z_1}-(cz_2^2+ab z_1z_2+z_2-z_1)\frac{\partial}{\partial z_2}\label{errz}
\end{equation}
is no rationally integrable. In fact, the field~\eqref{errz} is one of the equivalent fields defining the same foliation.

\item In the case $a_1 = a_0b_0/3$ (equivalently, $c=-ab/3$), the rational solution of the Risch equation is given by
${\ds y(x)= \frac{-\frac{6}{b}x^2 + 2}{x^3} }$. Thus, hypothesis (H2) is not fulfilled and the method does not decide.

\end{itemize}

%\end{example}

\bigskip

\noindent\textbf{Example 2 (continuation).}
%\begin{example}
%\label{ex:infin}
Take $n=m=0$, $A=a$ and $B=b$ in equation~\eqref{eq:infinity}. Due to $n<k-1$, we fall in case 1. Thus, $\widetilde{h}=k$ and $Y(x)=y_kx^k+\ldots+y_1x+2$. If condition \eqref{eqka4} is satisfied, then $2k-2\neq k$, which implies that we have not rational solution for all $k>2$.

Now, since $k=2$ does not satisfy condition~\eqref{eqka4}, we seek for $Y(x)$ using equations~\eqref{eqka} or~\eqref{eqka3}. We see that $\widetilde{h}=2$ and $Y(x)=y_2x^2+y_1x+2$. Now, by equation~\eqref{eqka3}, we obtain $$y_1=\frac{4}{a} ,\quad y_2=\frac{2ab+4}{a^2}.$$ 
The solution is given by 
\[
y(x) =\frac{2}{x^2}+{\frac {4}{ax}}+{\frac {2b}{a}}+{\frac{4}{a^2}}+{\rme}^{\frac {a}{x}}{c_1}.
\]
Thus, for this infinite family, the Risch equation has a rational solution for $k=2$ and no rational solution for $k>2$. 
Therefore, it follows that the vector field
\[
\widetilde{X}=\widetilde{P}(z_1,z_2) \frac{\partial }{\partial z_1}+\widetilde{Q}(z_1,z_2) \frac{\partial }{\partial z_2},
\]
with
\begin{equation*}
\widetilde{P}(z_1,z_2)=\sum_{i=0}^{N_1} P_i(z_1,z_2), \qquad
\widetilde{Q}(z_1,z_2)=\sum_{i=0}^{N_2}  Q_i(z_1,z_2),
\end{equation*}
where $P_i(z_1,z_2)$, $Q_i(z_1,z_2)$ are homogeneous polynomials of degree $i$, and satisfying conditions~\eqref{prop:Qfamily} is not rationally integrable if $P_N(z_1,z_2)=az_1^N$  and $P_{N-1}(z_1,z_2)=bz_1^{N-1}z_2-z_1^{N-1}$.

\bigskip

\noindent\textbf{Example 3 (continuation).}
Take $k=2$, $A=a_1x+a_0$ (with $a_1\neq 0$), $B\equiv 0$, $n=1$, and 
\[
\alpha=\frac{a_1x+a_0}{x^2},\qquad \beta=\frac{2a_1x+2a_0}{x^4},
\]
in equation~\eqref{whitf}. Since $n=k-1$ (case 2
of the algorithm) and $m<1$ we should consider two options. The first
one is when $a_1-2$ is not a negative integer, that is $a_1\in
[2,\infty)$, we fall in subcase 2.a and therefore the
Risch differential equation has no rational solutions.

Now, the second option is when $a_1-2$ is a negative integer, that is
$a_1\in (-\infty,1]\cap \mathbb{Z}$ and we fall in subcase 2.d.
Thus, $\widetilde{h}=2-1>0$ and then the Risch differential equation has always
one rational solution because $y_m,\ldots,y_1$ satisfy the algebraic
equation \eqref{eqka3} and we can obtain explicitly the polynomial $Y(x)$, which lead us to the rational solution $y(x)=\frac{Y(x)}{x^2}$.
For instance, setting $a_1=1$ and $a_0=a$, we obtain
$\widetilde{h}=1$ and $Y(x)=\frac{4}{a}x+2$ and the rational solution
of the Risch differential equation is
\[
y(x)=\frac{4}{ax}+\frac{2}{x^2}.
\]
Coming back to the variables $(z_1,z_2)=(1/y,x/y)$ the field~\eqref{whit} is (equivalent to)

\begin {equation}\widetilde{X}=(a_0 z_1^2+a_1z_1z_2)\frac{\partial}{\partial z_1}-(a_1z_2^2+a_0 z_1z_2-z_2^2+z_1)\frac{\partial}{\partial z_2}.\label{whittaz}\end{equation}
Thus we proved that for $a_1\notin (-\infty,1]\cap \mathbb{Z}$ the field
\eqref{whittaz} is not rationally integrable.
We observe that for $a_1=1$ it reduces to a field equivalent the following linear one
\begin{equation} 
\widetilde{X}=(a_0 z_1+z_2)\frac{\partial}{\partial z_1}-(a_0 z_2+1)\frac{\partial}{\partial z_2}.
\label{line}
\end{equation}

%\end{example}

\bigskip

\section{Conjecture: Extensions  to other type of first integrals}
\label{sec:conjecture}

We believe that our results provide obstructions for a more general kind of first integrals. More precisely, it seems natural to think assumptions (H1) and (mainly) (Hk) ($k\geq 2$) are not compatible with the existence of an {\it elementary} first integral (for a precise definition of it, see for instance~\cite{PS}). Indeed, the integral $\theta_k$ is too ``transcendent" to allow an elementary first integral.
Then  we state the following

\medskip
\noindent {\bf CONJECTURE}. {\it Under assumptions (H1) and (Hk) (for some $k\geq 2$)  the polynomial field has not an elementary first integral}.

\medskip

\noindent We illustrate this conjecture with the field \eqref{whit} with $a_1=1$. Since equation~\eqref{line} is linear in the variables $z_1,z_2$,  it must  be  integrable in some reasonable  sense.

\medskip

\noindent\textbf{Example.}
The function
\[
H(x,y)=\frac{ax+y}{y}\ {\rme}^{-\frac{a(x+a)}{ax+y}}
\]
is a first integral of the field
\[
X=(x^2-y) \frac{\partial}{\partial x}+ y(x+a)\frac{\partial}{\partial y}.
\]
Notice that $H$ is elementary, in fact of Darboux type (see~\cite{Da,DLLA}).

From the previous section we know that assumption (H1) is satisfied but not (H2), i.e,  $\theta_1$ is elementary. We compute $\theta_{k-1}$ 
up to  $k=10$ and all them are elementary (we skip the details), i.e, no assumption (Hk) is satisfied for $k\leq 10$ and hence until $k=10$ no obstruction to rational integrability is obtained, despite of the fact that the first integral is not rational but elementary.

\bigskip

\section*{Appendix A: Some \textsc{Maple} computations}

We will give some flavour of the analytical approach using \textsc{Maple} on two of the examples of Section~\ref{sec: exa}: Example~1 and Example~3. The results are in agreement with the algebraic results of Section~\ref{sec:alg}.

\bigskip

\noindent\textbf{Example 1.}
%\label{ex:errorf}
There are two possible cases: (i) $c\neq -ab/3$ and (ii) $c=-ab/3$.
\begin{itemize}

\item[(i)] If $b\neq 0$ and $c\neq -ab/3$ the solution of $\LVE_2$ is
\[
\begin{array}{cl}
\chi_1(x)&=c_1e^{\frac{2cx+b}{x^2}} x^2,\\
\chi_2(x)&=c_2  e^{\frac{2cx+b}{2x^2}} x+c_1\sqrt{2\pi}\ \frac{{\rm erf}\left(\frac{\sqrt{2}(cx+b)}{2\sqrt{-b}x}\right)(ab+3c)
e^{\frac{-c^2x^2+2bcx+b^2}{2bx^2}}x+\sqrt{2}\sqrt{-b}e^{\frac{2cx+b}{x^2}}(3x^2-b)}{(-b)^{3/2}x},
\end{array}
\]
and notice that for $b\neq 0$ the \emph{error function} ${\rm erf}\left(\frac{\sqrt{2}(cx+b)}{2x\sqrt{-b}}\right)$ is not
an elementary function (see, for instance,~\cite[page 48]{RITT}).
We remark that it is also possible  to express the error function by means of the Whittaker function (\cite{WHI}) $W_{-1/4,1/4}$ (and some elementary functions). Then by proving
that this Whittaker  function is not elementary, we obtain an alternative proof of the non-elementary character of the error function.
Notice that $G_2=Gal(\LVE_2)$ is connected and non abelian. Therefore, the vector field $X$ is
no rationally integrable.

\item[(ii)] In the case $c=-ab/3$, the solution of $\LVE_2$ is
\begin{eqnarray*}
\chi_1(x) &=& c_1 \,{x}^{2}{\rm e}^{-{\frac {b \left( 2\,ax-3 \right) }{3{x}^{2}}}}, \\
\chi_2(x) &=&
\frac{1}{bx} 
\left( 2 c_1 {\rme}^{-\frac{b (2ax-3)}{6x^2}} (-3x^2+b) + c_2 b x^2 \right)
{\rme}^{-\frac{b(2ax-3)}{6x^2}}.
\end{eqnarray*}
Solutions are elementary and its associated Risch equation has a rational solution.

\end{itemize}

\bigskip

\noindent\textbf{Example 3.}
%\label{ex:whit}
In this family the $\LVE_2$ becomes
\[
\begin{array}{cl}
\chi_1'&=2\dfrac{a_1x+a_0}{x^2}\chi_1\\[1.2ex]
\chi_2'&=\dfrac{a_1x+a_0}{x^2}\chi_2+2\dfrac{a_1x+a_0}{x^{4}}\chi_1.
\end{array}
\]
\textsc{Maple} gives the following  solution of $\LVE_2$
\[
\begin{array}{rcl}
\chi_1(x)&=& c_1x^{2a_1}e^{-2\frac{a_0}{x}}\\
\chi_2(x)&=& c_2 x^{a_1}
e^{-\frac{a_0}{x}}+c_1\dfrac{2a_0^2\ e^{\frac{-2a_0}{x}}x^{2a_1-2}+4a_0\ e^{\frac{-2a_0}{x}}x^{2a_1-1}
-4\left({\frac{a_0}{x}}\right)^{\frac{a_1}{2}}M\ e^{-\frac{3a_0}{2x}}x^{2a_1}}{a_0^2},
\end{array}
\]
where we have denoted by
\[
M:=M_{\kappa,\mu}(z), \quad \mbox{with} \quad \kappa=-\frac{a_1}{2}, \
 \mu=\frac{1-a_1}{2}, \ z=\frac{a_0}{x},
\]
the M Whittaker function solution of the Whittaker equation
\[
 y''-\left(\frac{1}{4}-\frac{\kappa}{z}+\frac{4\mu^2-1}{4z^2}\right)y=0.
\]
In all the checked values for $a_1\in [0,1]\cap Z$ the M function is elementary and it is not elementary for the rest of values, in agreement with our previous results.

\bigskip

\section*{Appendix B: The Galoisian Correspondence}

One of the key theorems of the Galois theory of linear differential equations, the Picard--Vessiot theory, is the existence of a Galoisian correspondence between intermediate differential field in the Picard--Vessiot extension of a linear differential system and algebraic subgroups of the Galois group.   It means that the structure of the solutions of the differential equation (i.e., how to solve the differential equation)  is a reflection of the  structure of Galois group.   An analogous theorem is valid for the classical Galois theory of polynomials and as there,  this theorem plays an essential r\^ole not only in the theory but also in  applications. For this reason this theorem is also called the Fundamental Theorem of Picard--Vessiot theory  (some references are~\cite{marram, CH,SV}).

We need some terminology and notations. Let $K\subset L$ be an  extension of differential fields. Then:
\begin{itemize}
\item  $Gal(L/K)$ is the group of differential automorphisms  of $L$ which are the identity over $K$.
\item If $H$ is a subgroup of $Gal(L/K)$, the intermediate field $K\subset L^H\subset L$ fixed by $H$ is $L^H:=\{ a\in L:\, \sigma(a)=a,\, \forall \sigma\in H\}$ .
%\item \textcolor{blue}{The extension  $K\subset L$ is normal if for any $a\in L\setminus K$ there exists  some $\sigma\in Gal(L/K)$, such that $\sigma(a)\neq a$ or, equivalently, if any element $a\in L$ fixed by  the complete group $Gal(L/K)$ necessary remains in $K$.}
\item Given a linear differential system
\begin{equation}
\frac {d{\bf y}}{dx}=A{\bf y}, \label{lin}
\end{equation}
with coefficients on a differential field $K$ (that is $A\in {\text Mat}(n,K)$), a \emph{Picard--Vessiot extension} $K\subset L$ of equation \eqref{lin} is the extension  generated by the  elements of a fundamental matrix $\Phi (x)=(\phi_{ij}(x))$,  $L=K(\phi_{11},...,\phi_{nn})$. From equation~\eqref{lin},  the field $L$ is a differential field.
\end{itemize}

%We recall that the Picard--Vessiot extensions are normal and that its Galois group $G=Gal(L/K)$ is a linear  algebraic group over the constants that along this paper becomes  the complex numbers.

\noindent\textbf{Theorem} (Correspondence in Picard-Vessiot theory)

\textit{
Let $K\subset L$ a Picard--Vessiot extension of \eqref{lin} and $G=Gal(L/K)$. Then there exists a one-to-one correspondence between the intermediary differential fields $K\subset S\subset L$ and the algebraic subgroups $H\subset G$ such that:
\begin{itemize}
\item[1)] Given $S$, $H=Gal(L/S)$.
\item[2)] Given $H$, $L^H=S$.
\item[3)] Picard--Vessiot  extensions  $K\subset S$ correspond to normal subgroups $H\subset G$ and then $G/H\approx Gal(S/K).$
\item[4)] If $M$ is a subgroup $M\subset G$, then to the differential field $L^M$ it corresponds the Zariski adherence of $H$, i.e., $Gal(L/L^M)=\bar{H}.$
\end{itemize}
}

\bigskip

\section*{Appendix C: the Galois grup of the third order variational equation \label{gal3}}

\subsection*{Linearising the variational equations}

Variational equations and their linearised expression play a key r\^ole in the theory of integrability of differential equations and, in particular, they have been also crucial in the results of this work. However they are very well known and commonly used in many textbooks on ODEs, higher order variational equations and (a possible, there can be several)
their linearisation are not so frequent. The aim of this annex  is to remind them in order to do this paper a bit more self-contained. We illustrate the method with the third order variational equation.

To start, let us assume that $y=\varphi(x)$ is a solution of $y'=f(x,y)$, $y(x_0)=\varphi(x_0)=y_0$. One possible way to derive the variational equations related to the $y$-variable (that is, \emph{freezing} $x_0$ as a initial $x$-value and moving transversally in the $y$-direction, that is, $y(x_0)=y_0+\eps$) is to look for solutions of such equation in the perturbative form, more suitable for computations  that the way we deduced the $\VE_2$ in the introduction,
\[
y(x)=\vfi(x) +  \eps \vfi_1(x) + \frac{\eps^2}{2!} \vfi_2(x) + \cdots + \frac{\eps^j}{j!} \vfi_j(x) + \cdots
\]
Substituting it into $y'=f(x,y)$ and expanding in Taylor series we get
\begin{eqnarray*}
\lefteqn{\vfi'(x) + \eps \vfi_1'(x) + \frac{\eps^2}{2!} \vfi_2'(x) + \frac{\eps^3}{3!} \vfi_3'(x) + \cdots = }\\
&& f(x,\vfi) + f_y(x,\vfi) \left( \eps \vfi_1 + \frac{\eps^2}{2!} \vfi_2 + \frac{\eps^3}{3!} \vfi_3 + \cdots \right) + \\
&& \frac{1}{2!} f_{yy}(x,\vfi) \left( \eps \vfi_1 + \frac{\eps^2}{2!} \vfi_2 + \frac{\eps^3}{3!} \vfi_3 + \cdots \right)^2 + \\
&& \frac{1}{3!} f_{yyy}(x,\vfi) \left( \eps \vfi_1 + \frac{\eps^2}{2!} \vfi_2 + \frac{\eps^3}{3!} \vfi_3 + \cdots \right)^3 + \cdots
\end{eqnarray*}
where $f_y(x,y)=\partial f / \partial y$, $f_{yy}=\partial^2 f / \partial y^2$, etc. Equating terms of power $\eps^k$ one obtains the variational equations of any order $k \geq 1$. Indeed, for $k=0$ one get the known solution $\vfi'(x)=f(x,\vfi(x))$. Concernng the rest of orders one has:
\[
\begin{array}{cl}
 & \vfi_1'     = f_y(x,\vfi) \vfi_1 \\
        & \vfi_1(x_0) = 1,                 \\\\
 & \vfi_2'     = f_y(x,\vfi) \vfi_2 + f_{yy}(x,\vfi) \vfi_1^2 \\
        & \vfi_2(x_0) = 0, \\\\
 & \vfi_3'     = f_y(x,\vfi) \vfi_3 + 3 f_{yy}(x,\vfi) \vfi_1 \vfi_2 + f_{yyy}(x,\vfi) \vfi_1^3 \\
        & \vfi_3(x_0) = 0,
\end{array}
\]
and so on. The initial values come from the condition $y(x_0)=y_0 + \eps$.
Observe that its solution exhibits a triangular scheme: we substitute the functions $\vfi_j$, $j=1,\ldots,k-1$ previously obtained in the equation for $\vfi_k$ and solve for it. Thus, all of them are non-homogeneous linear differential equations except the first one, $\VE_1$, which is homogeneous. They can be solved recurrently using the exponential of a primitive (for the first variational equation) and primitives.

Nevertheless they can be computed (specially, numerically) in this form, their expression is not the most suitable when regarding them in terms of differential Galois theory. In that framework, the common approach is to introduce intermediate variables which lead equations $\VE_k$, $k\geq 2$, into linear homogeneous systems. Namely, for the $\VE_3$ let us define the following variables
\[
 \chi_1=\vfi_1^3, \quad \chi_2 = \vfi_1 \vfi_2, \quad
\chi_3 = \vfi_3.
\]
Thus, the third order variational equation $\VE_3$ becomes  the $\LVE_3$
\begin{eqnarray*}
\chi_1'&=& 3 \vfi_1^2 \vfi_1'= 3 f_y(x,\vfi) \vfi_1^3 = 3 f_y(x,\vfi) \chi_1 \\
\chi_2' &=& \vfi_1'\vfi_2 + \vfi_1 \vfi_2' = f_y(x,\vfi) \vfi_1 \vfi_2 + \vfi_1 \left( f_y(x,\vfi) \vfi_2 + f_{yy}(x,\vfi) \vfi_1^2 \right) \\
        &=& 2 f_y(x,\vfi) \chi_2 + f_{yy}(x,\vfi) \chi_1 \\
\chi_3'&=& f_y(x,\vfi) \chi_3 +  3 f_{yy}(x,\vfi) \chi_2+f_{yyy}(x,\vfi) \chi_1 , \\
\end{eqnarray*}
which gives rise to the following triangular system
\[
 \hspace{1cm}
\left(
\begin{array}{c}
\chi_1' \\
\chi_2' \\
\chi_3'
\end{array}
\right) = \left(
\begin{array}{ccc}
3 f_y(x,\vfi)   &   0             &     0  \\
f_{yy}(x,\vfi)  & 2f_y(x,\vfi)    &     0  \\
f_{yyy}(x,\vfi) & 3f_{yy}(x,\vfi) &  f_y(x,\vfi)
\end{array}
\right)
\left(
\begin{array}{c}
\chi_1 \\
\chi_2 \\
\chi_3
\end{array}
\right).
\]
For variational equations of higher order a similar procedure can be performed.

\subsection*{Galois group of the third order linearised variational equation \boldmath{$\LVE_3$}}

To ease the notation in the computation of the Galois group associated to the third linearised variational equation $\LVE_3$, we introduce the following functions:
\begin{equation}
\alpha(x)=f_y(x,\vfi(x)), \qquad  \beta_2(x)=f_{yy}(x,\vfi(x)), \qquad \beta_3(x)=f_{yyy}(x,\vfi(x)),
\end{equation}
and define
\begin{equation*}
\omega={\rme}^{\int \alpha \,dx}, \qquad
\theta_1=\int \beta_2 {\rme}^{\int \alpha \,dx} \, dx, \qquad
\theta_2=\int \beta_3 {\rme}^{2\int \alpha \,dx} \, dx.
\end{equation*}
Observe, from such definitions, that we can just write
\begin{equation}
\label{def:omega:theta}
\omega={\rme}^{\int \alpha \,dx}, \qquad
\theta_1 =  \int \beta_2 \omega, \qquad \theta_2 = \int \beta_3 \omega^2.
\end{equation}
To solve the linear differential equations  $\LVE_3$ (homogeneous the first one and nonhomogeneous the others) we use that the general solution of a linear ODE $y'=a(x) y + b(x)$ is given, using the formula of variation of parameters, by ${\ds y=c {\rme}^{\int a} + {\rme}^{\int a} \int b {\rme}^{-\int a} }$. Indeed, we have
\begin{eqnarray*}
&& \chi_1 = c_1 \omega^3,\\
&& \chi_2 = c_2 \omega^2 + c_1 \omega^2 \theta_1, \\
&& \chi_3 = c_3 \omega +  3 c_2 \omega \theta_1+c_1\omega \theta_2 +\frac{3}{2} c_1 \omega \theta_1^2,
\end{eqnarray*}
where in the last equation it has been used that
\[
\int \beta_2 \omega \theta_1 =
\int \beta_2 {\rme}^{\int \alpha} \int \beta_2 {\rme}^{\int \alpha} = \int \theta_1 \theta_1' = \frac{1}{2} \theta_1^2= \frac{1}{2} \left( \int \beta_2 {\rme}^{\int \alpha} \right)^2
\]
and with $c_1,c_2,c_3$ arbitrary complex constants.
That is, a possible fundamental matrix of system $\LVE_3$ is
\[
\Phi_3=\left(
\begin{array}{ccc}
\omega^3 &  0 & 0\\
\omega^2 \theta_1 & \omega^2 & 0 \\
\frac{3}{2}\omega\theta_1^2+ \omega \theta_2  & 3 \omega \theta_1 & \omega
\end{array}
\right).
\]
For $\sigma \in G_3$, the Galois group of $\LVE_3$ we have $\sigma(\omega)=\lambda \omega$, $\sigma(\theta_1)=\lambda \theta_1+\mu$,
$\sigma(\theta_2)=\lambda^2 \theta_2+{\nu}$, with $\lambda\neq 0$, $\mu,\nu$ suitable complex constants. It follows that the corresponding Galois group $G_3$ is represented by an algebraic subgroup of the  triangular group given by the matrices
is contained in
\[
B_\sigma=\left(
\begin{array}{ccc}
\lambda^3 & 0 & 0 \\
\lambda^2\mu & \lambda^2 & 0 \\
\lambda\nu+\frac{3}{2}\lambda\mu^2 & 3\lambda\mu & \lambda
\end{array}
\right),
\]
analogously as it was done for the linearised second order variational equation $\LVE_2$.
Therefore, assumption (H1) is not verified if and only if $\lambda=1$ and assumptions (H2) and (H3) are not verified if and only if $\mu=0$ and $\nu=0$, respectively.

\section*{Acknowledgements}

\noindent JM is member of the Universidad Polit\'ecnica de Madrid research group ``Modelos
Matem\'aticos no lineales''. His work has been partially supported by this research group.
CP has been partially supported by Spanish MECC-FEDER grant MTM2015-65715-P and MINECO-FEDER grant MTM2016-77278-P, and by the Catalan AGAUR grant 2014SGR-504. JTL has been supported by Spanish MECC-FEDER grant MTM2015-65715-P, the Catalan AGAUR grant 2014SGR-504 and the grant 14-41-00044 of RSF at the Lobachevsky University of Nizhny Novgorod (Russia). PB thanks Universidad Sim\'on Bol\'ivar (Barranquilla, Colombia) for its support during the final stage of this work. All the authors thank the EPSEB (UPC) for its support.

And last, but not least, the authors would like to thank the people from the UPM Integrability Seminar in Madrid, the UB-UPC Dynamical Systems Seminar in Barcelona and the Dynamical Systems Seminar at the UAB in Bellaterra, where preliminary versions of this work
were exposed. In particular to C. Sim\'o and A. Delshams for their interesting remarks and suggestions.

\end{document}